# A GPU-accelerated Cartesian grid method for the heat, wave and Schrödinger equations on irregular domains

Liwei Tan[1], Minsheng Huang[1], Wenjun Ying[1,2,*]

[1] *School of Mathematical Sciences, Shanghai Jiao Tong University, Shanghai 200240, P. R. China.*

[2] *School of Mathematical Sciences, MOE-LSC and Institute of Natural Sciences, Shanghai Jiao Tong University, Minhang, Shanghai 200240, P. R. China.*

**Abstract.** Based on Ying's kernel-free boundary integral(KFBI) method[1], a second-order method for general elliptic partial differential equations (PDEs), this paper develops a GPU-accelerated KFBI method for the heat, wave and Schrödinger equations on the irregular domain. Since the limitation of time steps imposed by CFL conditions in the explicit scheme and the inadequate accuracy generated by the fully implicit scheme for the Laplacian operator, the paper selects a series of second-order time discrete schemes, and the Laplacian operator is split into explicit and implicit mixed ones. The Crank-Nicolson method is used to discretize the heat equation in temporal dimension while the implicit $\theta-$scheme is for the wave equation. The Strang splitting method is applied to the Schrödinger equation. After discretizing the temporal dimension implicitly, the heat, wave and Schrödinger equations are transformed into a sequence of elliptic equations. The Laplacian operator on the right-hand side of the elliptic equation is obtained from the numerical scheme instead of being discretized and corrected by the five-point difference method. A Cartesian grid-based KFBI method is used to solve the resulting elliptic equations. The KFBI method is accelerated by the graphics processing unit (GPU) with a parallel Cartesian grid solver, achieving a high degree of parallelism. Numerical results show that the proposed method has a second-order accuracy for the heat, wave, and Schrödinger equations. Additionally, the GPU-accelerated solvers for the three types of time-dependent equations are 30 times faster than CPU-based solvers.

# 1 Introduction

The heat, wave, and Schrödinger equations encountered many scientific and engineering applications. Examples include biological systems and option pricing associated with the heat equation [2–5], electromagnetic wave propagation and underwater acoustics related to the wave equations [6–8], and quantum mechanics issues associated with the Schrödinger equation [3,4,9,10]. Over the past few decades, there have been extensive studies in numerical methods for those three types of time-dependent PDEs, including the finite element method [11,12], the finite difference method [13,14], the pseudo-spectral method [15] and the boundary integral method(BIM) [16,17]. The primary purpose of this work is to present an efficient GPUaccelerated and second-order accurate kernel-free boundary integral (KFBI) method combined with an implicit time discretization scheme as an alternative for the approaches above in solving the heat, wave and Schrödinger equations.

* Corresponding author, Wenjun Ying: wying@sjtu.edu.cn

One can distinguish three approaches of BIM to the application of the heat, wave and Schrödinger problems, space-time boundary integral equation methods [17–20], Laplace-transform methods [17–21], and time-stepping methods [16,17,20,22,23]. For the space-time boundary integral method, the original equation is reformed as a space-time boundary integral equation by Green's third identity, then the boundary element method(BEM) [24], Galerkin method [25] and collocation method [26] can be applied to solve the space-time boundary integral equation. Remarkably, because the numerical methods constructed from these space-time boundary inte gral equations are global in time, the solution can be obtained in one step for the entire time interval [20]. However, the difficulty of obtaining Green's functions for general domains and the high cost of solving the large discrete system matrix and higher-dimensional integrals prevent it from becoming a universal method [20]. For the Laplace transform method, the primal problem is transformed to an elliptic boundary value problems in the frequency domain with an eigen value parameter λ depending on the frequency ω by taking the Laplace transform of the variable. Then the elliptic boundary value problems in the finite frequency domain are solved by the BIM. Finally, the solution in the frequency domain space is converted back to the time domain space through the inverse Laplace transform. Due to the good characteristics of the Laplacian operator, this method's solution process is relatively simple [27]. However, there will be some ωl with large absolute value, which means large negative real parts for $\lambda(\omega_l)$. An alternative to solving the three types of time-dependent equations is first to discretize the temporal dimension by an implicit scheme, then solve the resulting elliptic equations for each time step by boundary integral method. Significantly, a careful treatment of time discretization and spatial discretization at each time step is of utmost importance to ensure the long-term numerical stability of the discrete scheme.

Since the hyperbolicity of wave and Schrödinger equations, it is expected that the discrete schemes employed in time-stepping methods respect the dissipative properties enjoyed by the original systems. Structure-preserving algorithms are achieved by constructing numerical methods, which can preserve some properties of continuous systems [28] and designing numerical schemes that preserve dissipative properties are crucial for accurate and reliable long-time numerical simulations of dissipative systems. For the wave equation, the conventional methods, such as Lax-Friedrichs and Lax-Wendroff methods, are nonsymplectic, which causes oscillation [29]. Dufort-Frankel scheme and Saulyev scheme are appropriate but they could produce a dissipative phenomenon[30]. An alternative method to solve the wave equation is the implicit $\theta-$scheme. When $\theta = 1/4$, the unconditionally stable scheme is derived from the Euler midpoint scheme, which is a simple symplectic scheme. The leap-frog method is non-dissipative. but there are potential stability issues if it is used for variable coefficients or nonlinear problems. For the Schrödinger equation, the structure-preserving splitting methods, for instance, the first-order Godunov splitting method [32], the second-order Strang splitting method [33] and the fourth-order Yoshida method [34] are three common methods. For the structure-preserving splitting methods, the time-splitting sine pseudo-spectral(TSSD) method [35,36] and the timesplitting finite difference (TSFD) method [37] are unconditionally stable, time-reversible, and time-transverse invariant. Nevertheless, for non-smooth solutions, the spectral method could result in a loss of precision. For the finite difference time domain(FDTD) method [38], the CrankNicolson finite difference(CNFD) method[39], which is unconditionally stable, conserves the mass and energy in the discretized level. However, more significant dispersion errors in a longtime simulation are the bottleneck yet to be overcome.

The main difficulty in the time-stepping method is solving a sequence of elliptic equations efficiently. The BIM and the BEM [40–45] offer several notable advantages when dealing with elliptic problems on irregular domains [5,46,47]. The reduced dimensionality of the model and the ease of complex boundary capturing are two evident merits. Nonetheless, some obvious drawbacks, such as the difficulty in obtaining Green's functions for general domains, the complexity of the mathematical analysis and error analysis [48], which restrict the conventional BIM from being a general-purpose method. To avoid directly representing Green's function, Mayo [49–51] evaluates a boundary or volume integral by solving an equivalent but simple interface problem on a Cartesian grid with FFT-based solvers. Inspired by Mayo's work, in 2007, Ying and Henriquez [1] proposed the KFBI method for elliptic boundary value problems as a competitive alternative to other Cartesian grid methods such as the high-order matched interface and boundary method [52], the ghost fluid method [53], the immersed boundary(IB) method [54], the volume of fluid(VOF) method [55], the decomposed II method [56] and so on [57–60]. Furthermore, Ying and Wang [61] applied their method to solve variable coefficients elliptic PDEs in 2014. In brief, the KFBI method has a few advantages[1]. Firstly, The KFBI method can be applied to smooth irregular domains and works with more general elliptic operators with possible anisotropy and inhomogeneity. Secondly, The KFBI method can easily work with Cartesian grids so that fast elliptic solvers such as the FFT-based solvers are applicable. In addition, Ying and Henriquez [1] show second-order convergence rates with smaller local truncation errors on uniform grids, while it has fourth-order accuracy presented by Xie and Ying [62]. Last but not least, it preserves the symmetry and positive definiteness of the coefficient matrix of the discrete system [1]. However, due to the hyperbolicity of the wave equation and the Schrödinger equation,the KFBI method is not generalized to these two types of equations yet. In addition, because the computational work is essentially linearly proportional to the number of nodes on the Cartesian grid, it costs a lot of CPU time when performing large-scale 2D and 3D computations [62,63].

The Cartesian grid methods are well-suited for parallel computing on GPUs. Firstly, Cartesian grids have a straightforward advantage in data representation, where each grid node is typically assigned to a dedicated thread for computation, resulting in high parallelism. Secondly, Cartesian grids possess superficial topological relationships, often relying solely on self and neighboring point information during the computation process. This characteristic is conducive to efficient data transfer, communication, and memory sharing. Several works have addressed the GPU acceleration of Cartesian grid methods in recent years [64–67]: the GPU-accelerated VOF by Rajesh Reddy and R. Banerjee [64], the CUDA-Based IB method by S. K. Layton, A.Krishnan and L. A. Barba [65], the TVD Runge–Kutta method on multiple GPUs by Liang. S, Liu. W and Yuan. L [66], the multi-GPU-based lattice Boltzmann algorithm by Huang. C, Shi. B, He. N and Chai. Z [67]. The KFBI method is also a Cartesian grid method, promising to design a GPU-accelerated version for significantly enhancing computational efficiency. In fact, due to the simple mesh topology in the Cartesian grid, creating a GPU-accelerated Cartesian grid solver with a high degree of parallelism is not tricky.

This paper applies the KFBI method combined with the implicit time discretization scheme to solve the heat, wave and Schrödinger equations. Since the second-order KFBI method is applied in space directions, in this paper the second-order symplectic time integrator is selected to keep the time order consistent with the space order. The Crank-Nicolson method, the implicit theta-scheme and the Strang splitting scheme are respectively employed to solve the heat wave and Schrödinger equations. First, the heat, wave and Schrödinger equations are discretized in temporal dimension implicitly, and then an elliptic boundary value problem on the irregular domain is solved in each time step. The right-hand side of the elliptic equation is considered a new variable and directly evaluated by the transformed time-advancing scheme. With the KFBI method, evaluation of a boundary or volume integral is replaced by interpolation of a Cartesian grid-based solution, which satisfies an equivalent interface problem. To solve the interface problem efficiently, the KFBI method is combined with the parallel technology in this paper. The significant contributions of our work can be concluded:(1)the KFBI method is applied to solve the three types of PDEs in 2D and 3D. (2)the KFBI method is implemented on a single GPU.

The remainder of this paper is organized as follows. Initial-boundary value problems of the heat equation, the wave equation and the Schrödinger equation are described in section 2. The numerical methods, including the implicit time discretization schemes and the GPU-accelerated KFBI method, are summarized in section 3. The numerical solutions are presented in section4. Finally, a short discussion is made in section 5.

# 2 The heat, wave and Schrödinger equations

Supposing $\Omega$ is a bounded two-dimensional irregular and complex domain whose boundary $\Gamma = \partial\Omega$ is at least twice continuously differentiable. Let $u(\mathbf{x}, t) = u(x, y, t)$ be an unknown function of $\mathbf{x} \in \mathbf{R}^2$. Assuming $u_0(\mathbf{x})$, $\partial_t u_0(\mathbf{x})$ are the initial data defined in $\Omega$, $g_D(\mathbf{x}, t)$ or $g_N(\mathbf{x}, t)$ are known functions of $\mathbf{x}$ with sufficient smoothness. $\partial_\mathbf{n} u(\mathbf{x}, t)$ denotes the normal derivative of $u(\mathbf{x}, t)$ on the boundary, where $\mathbf{n}$ denotes the unit outward normal on $\Gamma$. For each equation, we consider a model problem of each type.

- The heat equation

Dirichlet and Neumann initial boundary value problems(IBVPs) of the heat equation are considered as follows:

$$\partial_t u(\mathbf{x}, t) = \Delta u(\mathbf{x}, t) \text{ in } \Omega \times (0, T] \quad (2.1)$$

- The wave equation

Dirichlet and Neumann initial boundary value problems(IBVPs) of the wave equation are considered as follows:

$$\partial_{tt} u(\mathbf{x}, t) = \Delta u(\mathbf{x}, t) \text{ in } \Omega \times (0, T]$$

- The Schrödinger equation

Dirichlet and Neumann initial boundary value problems(IBVPs) of the Schrödinger equation are considered as follows. Let $v(\mathbf{x})$ be a given potential function, and $w$ is a non-negative real parameter

$$-i\partial_t u(\mathbf{x}, t) + \Delta u(\mathbf{x}, t) + v(\mathbf{x})u(\mathbf{x}, t) + w|u(\mathbf{x}, t)|^2 u(\mathbf{x}, t) = 0 \text{ in } \Omega \times (0, T] \quad (2.2)$$

- The initial and boundary conditions

The initial boundary condition can be represented as:

$$u(\mathbf{x}, 0) = u_0(\mathbf{x}), \quad \mathbf{x} \in \Omega$$

The initial boundary condition of the first time derivative for the wave equation reads:

$$\partial_t u(\mathbf{x}, 0) = \partial_t u_0(\mathbf{x}), \quad \mathbf{x} \in \Omega$$

The Dirichlet boundary condition can be considered as:

$$u(\mathbf{x}, t)|_\Gamma = \mathbf{g}_D(\mathbf{x}, t), \quad \text{for } t > 0$$

The Neumann boundary condition can be written as:

$$\partial_\mathbf{n} u(\mathbf{x}, t)|_\Gamma = \mathbf{g}_N(\mathbf{x}, t), \quad \text{for } t > 0$$

# 3 Numerical method

## 3.1 Semi-Discretization

The heat, wave and the Schrödinger equation are supposed to be discretized in temporal dimension first. The Crank-Nicolson method is used to discretize the heat equation, while the implicit $\theta$-scheme[68] is employed to discretize the wave equation. The Strang splitting method is applied to the Schrödinger equation. Let $u^n(\mathbf{x})$ be the approximate solution of $u(\mathbf{x}, t^n)$ and assume time step $\tau = t^{n+1} - t^n$. For the sake of simplicity, only the three types of equations with the Dirichlet boundary condition are shown below, while the Neumann condition is similar.

- Heat equation

It is well-known that the Crank-Nicolson scheme is an unconditionally stable implicit method with second-order accuracy. Applying the Crank-Nicolson scheme to the heat equation with the Dirichlet boundary condition (2.1), we obtain

$$\frac{u^{n+1}(\mathbf{x}) - u^n(\mathbf{x})}{\tau} = \frac{1}{2}(\Delta u^{n+1}(\mathbf{x}) + \Delta u^n(\mathbf{x})), \quad \mathbf{x} \in \Omega$$
$$u^0(\mathbf{x}) = u_0(\mathbf{x}), \quad \mathbf{x} \in \Omega \quad (3.1)$$
$$u^n(\mathbf{x}) = g_D(\mathbf{x}, t^n), \quad \mathbf{x} \in \Gamma$$

Rewriting the equation (3.1), the actual elliptic boundary value problem that one has to solve at each time step is therefore

$$\Delta u^{n+1}(\mathbf{x}) - \frac{2}{\tau} u^{n+1}(\mathbf{x}) = -\frac{2}{\tau} u^n(\mathbf{x}) - \Delta u^n(\mathbf{x}), \quad \mathbf{x} \in \Omega$$
$$u^n(\mathbf{x}) = g_D(\mathbf{x}, t^n), \quad \mathbf{x} \in \Gamma \quad (3.2)$$

Denoting $F^{n+1}(\mathbf{x}) = \frac{2}{\tau} u^n(\mathbf{x}) + \Delta u^n(\mathbf{x})$ and substituting it into the right-hand side of (3.2), we get

$$\Delta u^{n+1}(\mathbf{x}) = \frac{2}{\tau} u^{n+1}(\mathbf{x}) - F^{n+1}(\mathbf{x}) \quad (3.3)$$

It is empirical that there exists numerical instability when computing $u_{xx}^{n+1}, u_{yy}^{n+1}$. The right-hand-side of the elliptic equation is considered as a new variable and one can avoid computing those derivatives by substituting (3.3) into the next time step, which can be expressed as

$$F^{n+2}(\mathbf{x}) = -\Delta u^{n+2}(\mathbf{x}) + \frac{2}{\tau} u^{n+2}(\mathbf{x}) = \frac{4}{\tau} u^{n+1}(\mathbf{x}) - F^{n+1}(\mathbf{x}), \quad n = 0,1, \ldots \quad (3.4)$$

- Wave equation

The implicit $\theta$-scheme is an implicit scheme where $\theta$ is a parameter. Moreover, different $\theta$ leads to different schemes. A second-order, unconditionally stable scheme is derived when $\theta \in [1/4, 1/2]$. In detail, a general version of the implicit $\theta$-scheme of the wave equation can be expressed as

$$\frac{u^{n+1}(\mathbf{x}) - 2u^n(\mathbf{x}) + u^{n-1}(\mathbf{x})}{\tau^2} = \theta \Delta u^{n+1}(\mathbf{x}) + (1 - 2\theta)\Delta u^n(\mathbf{x}) + \theta \Delta u^{n-1}(\mathbf{x}), \quad \mathbf{x} \in \Omega \quad (3.5)$$
$$u^n(\mathbf{x}) = g_D(\mathbf{x}, t^n), \quad \mathbf{x} \in \Gamma$$

Setting the parameter $\theta = \frac{1}{4}$, and after simply manipulation of equation(3.5), it follows:

$$\Delta u^{n+1}(\mathbf{x}) - \frac{4}{\tau^2} u^{n+1}(\mathbf{x}) = -\frac{8}{\tau^2} u^n(\mathbf{x}) + \frac{4}{\tau^2} u^{n-1}(\mathbf{x}) - 2\Delta u^n(\mathbf{x}) - \Delta u^{n-1}(\mathbf{x}), \quad \mathbf{x} \in \Omega \quad (3.6)$$
$$u^n(\mathbf{x}) = g^n(\mathbf{x}), \quad \mathbf{x} \in \Gamma$$

Denoting $F^{n+1}(\mathbf{x}) = \frac{8}{\tau^2} u^n(\mathbf{x}) - \frac{4}{\tau^2} u^{n-1}(\mathbf{x}) + 2\Delta u^n(\mathbf{x}) + \Delta u^{n-1}(\mathbf{x})$ and substituting it into (3.6), we get

$$\Delta u^{n+1}(\mathbf{x}) = \frac{4}{\tau^2} u^{n+1}(\mathbf{x}) - F^{n+1}(\mathbf{x}) \quad (3.7)$$

Similar to the derivation of the equation(3.4), substituting (3.7) into (3.6) when $t = t^{n+2}$, we obtain

$$F^{n+2}(\mathbf{x}) = -\Delta u^{n+2}(\mathbf{x}) + \frac{4}{\tau^2} u^{n+2}(\mathbf{x}) = \frac{16}{\tau^2} u^{n+1}(\mathbf{x}) - 2F^{n+1}(\mathbf{x}) - F^n(\mathbf{x}), \quad n = 0, 1, \ldots$$

- The Schrödinger equation

Operator splitting is a numerical method for solving decomposable differential equations into a sum of differential operators. For the Schrödinger equation, this method is widely used because it allows us to effectively deal with the nonlocal Laplacian and the nonlinear terms simultaneously. As usual, we can write the Schrödinger equation in the following scheme

$$\partial_t u(\mathbf{x}, t) = A(u(\mathbf{x}, t))u(\mathbf{x}, t) + Bu(\mathbf{x}, t)$$

Where the linear operator $A = -i\Delta$ and the nonlinear part $B = -i(v(\mathbf{x}) + w|u(\mathbf{x})|^2)$. This abstract formulation will be convenient for clarifying the Strang splitting scheme below.

Also known as the symmetric operator splitting, the Strang splitting method can be written as

$$u(\mathbf{x}, t + \tau) = \mathcal{S}^A\left(\frac{\tau}{2}\right) \mathcal{S}^B(\tau) \mathcal{S}^A\left(\frac{\tau}{2}\right) u(\mathbf{x}, t)$$

Where $\mathcal{S}^A$ is the exact solution operators of the equation (3.8), and $\mathcal{S}^B$ is the exact solution operators of the nonlinear equation(3.9) respectively

$$iu_t(\mathbf{x}) = \frac{1}{2}\Delta u(\mathbf{x}) \quad (3.8)$$

$$iu_t(\mathbf{x}) = v(\mathbf{x})u(\mathbf{x}) + w|u(\mathbf{x})|^2 u(\mathbf{x}) \quad (3.9)$$

Now, we replace the exact solution operators $\mathcal{S}^A$ and $\mathcal{S}^B$ by their numerical approximations in space semi-discrete $(\mathcal{S}_h^A)^{forward}$, $(\mathcal{S}_h^A)^{backward}$ and $\mathcal{S}_h^B$.

**Step1:** $\mathcal{S}^A \to (\mathcal{S}_h^A)^{forward}$

For the linear equation (3.8), the forward Euler scheme can be written as

$$u^{n+1}(\mathbf{x}) = u^n(\mathbf{x}) - i\frac{\tau}{2}\Delta u^n(\mathbf{x})$$

**Step2:** $\mathcal{S}^B \to \mathcal{S}_h^B$

For the nonlinear equation (3.9), the temporal dimension is discretized implicitly, in particular, the Crank-Nicolson scheme can be expressed as

$$i\frac{u^{n+1}(\mathbf{x}) - u^n(\mathbf{x})}{\tau} = \frac{1}{2}(v(\mathbf{x})u^n(\mathbf{x}) + v(\mathbf{x})u^{n+1}(\mathbf{x})) + \frac{w}{2}(|u^{n+1}(\mathbf{x})|^2 u^{n+1}(\mathbf{x}) + |u^n(\mathbf{x})|^2 u^n(\mathbf{x}))$$

(3.10)

Simplify the equation (3.10) as

$$u^{n+1}(\mathbf{x}) + \frac{i\tau}{2}(v(\mathbf{x}) + w|u^{n+1}(\mathbf{x})|^2)u^{n+1}(\mathbf{x}) = u^n(\mathbf{x}) - \frac{i\tau}{2}(v(\mathbf{x}) + w|u^n(\mathbf{x})|^2)u^n(\mathbf{x})$$

The $u^{n+1}$ is regarded as an unknown variable in the equation, the nonlinear equation can be solved by the Newton method.

**Step3:** $\mathcal{S}^A \to (\mathcal{S}_h^A)^{backward}$

For the linear equation (3.8), the backward Euler method can be represented as

$$u^{n+1}(\mathbf{x}) = u^n(\mathbf{x}) - i\frac{\tau}{2}\Delta u^{n+1}(\mathbf{x}) \quad (3.11)$$

Rewrite the equation (3.11), the actual elliptic boundary value problem that one has to solve at each time step is therefore

$$u^{n+1}(\mathbf{x}) + i\frac{\tau}{2}\Delta u^{n+1}(\mathbf{x}) = u^n(\mathbf{x}) \quad (3.12)$$

Now, introducing two intermediate variables $u^*(\mathbf{x})$ and $u^{**}(\mathbf{x})$, the Strang splitting method is employed to solving the problem(2.2), which leads to the following algorithm:

$$\begin{cases} u^*(\mathbf{x}) = u^n(\mathbf{x}) - i\frac{\tau}{2}\Delta u^n(\mathbf{x}) \\ u^{**}(\mathbf{x}) + \frac{i\tau}{2}(v(\mathbf{x}) + w|u^{**}(\mathbf{x})|^2)u^{**}(\mathbf{x}) = u^*(\mathbf{x}) - \frac{i\tau}{2}(v(\mathbf{x}) + w|u^*(\mathbf{x})|^2)u^*(\mathbf{x}) \\ u^{n+1}(\mathbf{x}) + i\frac{\tau}{2}\Delta u^{n+1}(\mathbf{x}) = u^{**}(\mathbf{x}) \end{cases}$$

Note that the above algorithm can be rewritten in the abstract formula

$$u^{n+1}(\mathbf{x}) = (S_h^A)^{backward}\left(\frac{\tau}{2}\right) S_h^B(\tau) (S_h^A)^{forward}\left(\frac{\tau}{2}\right) u^n(\mathbf{x}),$$

## 3.2 Kernel-free boundary integral method

We first introduce the modified Helmholtz equation, which is of most importance in this work and has the same formulation as the (3.2),(3.6) and (3.12). Let $\mathbf{x} = (x, y) \in \mathbf{R}^2$ and $u(\mathbf{x}) = u(x, y)$ be an unknown function of $\mathbf{x} \in \mathbf{R}^2$. Let $\Omega \in \mathbf{R}^2$ be a bounded domain with at least twice continuously differentiable boundary $\Gamma = \partial\Omega$. The modified Helmholtz equation reads as follows

$$\Delta u(\mathbf{x}) - \kappa u(\mathbf{x}) = f(\mathbf{x}), \quad x \in \Omega \qquad (3.13)$$

subject to either Dirichlet boundary condition:

$$u(\mathbf{x}) = g_D(\mathbf{x}), \quad x \in \Gamma$$

or the Neumann boundary condition:

$$\partial_\mathbf{n} u(\mathbf{x}) = g_N(\mathbf{x}), \quad x \in \Gamma$$

Here, $f(\mathbf{x})$, $g_D(\mathbf{x})$ and $g_N(\mathbf{x})$ are known functions. $\kappa$ is assumed to be a positive constant for the modified Helmholtz equation in this paper by default; $\mathbf{n} = (n_x, n_y)$ is the unit outward normal on $\Gamma$; $\partial_\mathbf{n} u(\mathbf{x})$ denotes the normal derivative of $u(\mathbf{x})$ on the boundary $\Gamma$.

### 3.2.1 Boundary integral equation

As shown in Fig. 1, to solve the BVPs above by the KFBI method, we first embed the irregular domain $\Omega$ into a larger rectangle domain $\mathcal{B}$.

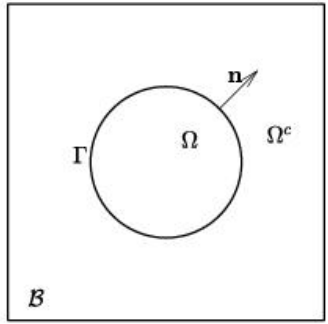

Figure 1: KFBI computation domain

According to the standard BIM[46, 48], let $G(\mathbf{x},\mathbf{y})$ be Green's function on the rectangle $\mathcal{B}$ associated with the elliptic PDE (3.13), which satisfies for $\mathbf{y} \in \mathcal{B}$,

$$\triangle G(\mathbf{x},\mathbf{y}) - \kappa G(\mathbf{x},\mathbf{y}) = \delta(\mathbf{x} - \mathbf{y}) \quad \mathbf{x} \in \mathcal{B}$$
$$G(\mathbf{x},\mathbf{y}) = 0 \quad \mathbf{x} \in \partial\mathcal{B}$$

Let $\mathbf{n_y}$ be the unit outward normal vector at point $\mathbf{y} \in \Gamma$, and $\varphi$ be the density function, define the double layer boundary integral by

$$(W\varphi)(\mathbf{x}) := \int_\Gamma \frac{\partial G(\mathbf{y},\mathbf{x})}{\partial \mathbf{n_y}} \varphi(\mathbf{y}) ds_\mathbf{y} \text{ for } \mathbf{x} \in \Omega \cup \Omega^c \quad (3.14)$$

Let $\psi$ be the density function, and define the single layer boundary integral by

$$(V\psi)(\mathbf{x}) := \int_\Gamma G(\mathbf{y},\mathbf{x}) \psi(\mathbf{y}) ds_\mathbf{y} \text{ for } \mathbf{x} \in \Omega \cup \Omega^c \quad (3.15)$$

As the right side $f(\mathbf{x})$ is only defined on $\Omega$, the volume integral is defined by

$$(Yf)(\mathbf{x}) := \int_\Omega G(\mathbf{y},\mathbf{x}) f(\mathbf{y}) d\mathbf{y} \text{ for } \mathbf{x} \in \mathbf{R}^2 \quad (3.16)$$

The Dirichlet BVPs can be reformulated as a Fredholm boundary integral equation of the second kind by Green's third identity

$$\frac{1}{2}\varphi(\mathbf{x}) + (W\varphi)(\mathbf{x}) + (Yf)(\mathbf{x}) = g_D(\mathbf{x}) \quad (3.17)$$

The solution $u(\mathbf{x})$ of the Dirichlet BVP can be reformulated by

$$u(\mathbf{x}) = (W\varphi)(\mathbf{x}) + (Yf)(\mathbf{x}) \quad (3.18)$$

The Neumann BVPs can be reformulated as a Fredholm boundary integral equation of the second kind[69, 70] by the Green's third identity

$$\frac{1}{2}\psi(\mathbf{x}) - \frac{\partial}{\partial \mathbf{n_x}}(V\psi)(\mathbf{x}) + \frac{\partial}{\partial \mathbf{n_x}}(Yf)(\mathbf{x}) = g_N(\mathbf{x}) \quad (3.19)$$

The solution $u(\mathbf{x})$ can be expressed as

$$u(\mathbf{x}) = (Yf)(\mathbf{x}) - (V\psi)(\mathbf{x}) \quad (3.20)$$

The BIE corresponding to the Dirichlet BVP (3.17) can be solved by the Richardson iteration:

$$\varphi_{k+1}(\mathbf{x}) = (1 - \frac{\gamma}{2})\varphi_k(\mathbf{x}) + \gamma[g_D(\mathbf{x}) - (Yf)(\mathbf{x}) - (W\varphi_k)(\mathbf{x})], \quad \mathbf{x} \in \Gamma \quad (3.21)$$

The BIE corresponding to the Neumann BVP (3.19) can be solved by the Richardson iteration:

$$\psi_{k+1}(\mathbf{x}) = (1 - \frac{\gamma}{2})\psi_k(\mathbf{x}) + \gamma[g_N(\mathbf{x}) - \frac{\partial}{\partial \mathbf{n_x}}(Yf)(\mathbf{x}) + \frac{\partial}{\partial \mathbf{n_x}}(V\psi)(\mathbf{x})], \quad \mathbf{x} \in \Gamma \quad (3.22)$$

For $k = 0,1,2, …$, respectively. Here $\gamma \in (0,1)$ is an iteration parameter. Choosing the appropriate initial value $\varphi_0$, $\psi_0$, the Richardson iteration converges within a fixed number of steps. The unknown density function $\varphi(\mathbf{x})$, $\psi(\mathbf{x})$ is obtained when the iteration (3.21) converges. After solving the linear system (3.21), the approximation of $u(\mathbf{x})$ can be calculated according to the formula (3.18), (3.20).

### 3.2.2 Equivalent interface problems

Due to the complexity of Green's function $G(\mathbf{x}, \mathbf{y})$ defined in the bounded domain $\mathcal{B}$, the main challenge in this process is on evaluation of boundary or volume integrals in (3.21), (3.22). To avoid directly discretizing the integrals in terms of Green's function by numerical quadratures, the boundary integral and volume integral are equivalent to solutions of the interface problems (3.24), (3.27), (3.29).

Assuming that any function $u(\mathbf{x})$ and its partial derivatives involved in the following are at least piece-wise smooth, with potential discontinuity only existing on the interface $\Gamma$. Denoting symbol "+" or "-" to distinguish the related functions restricted to the subdomain $\Omega$ or $\Omega^c$. $u^+$ ($u^-$) means the restrictions of $u$ in the subdomains $\Omega(\Omega^c)$

$$u^+(\mathbf{x}) = u(\mathbf{x})|_\Omega, \quad u^-(\mathbf{x}) = u(\mathbf{x})|_{\Omega^c}$$

For the point $\mathbf{x}$ near the interface, $u^+(\mathbf{x})$ and $u^-(\mathbf{x})$ represents the limit values of $u(\mathbf{x})$ from the corresponding side of the domain boundary

$$u^+(\mathbf{x}) = \lim_{z \to x, z \in \Omega} u(\mathbf{z}), \quad u^-(\mathbf{x}) = \lim_{z \to x, z \in \Omega^c} u(\mathbf{z})$$

For the partial derivatives, $u_x^+(\mathbf{x})$ ($u_y^+(\mathbf{x})$) and $u_x^-(\mathbf{x})$ ($u_y^-(\mathbf{x})$) are defined as the limit values of $u_x(\mathbf{x})(u_y(\mathbf{x}))$ from the corresponding side of the domain boundary

$$\begin{aligned} u_x^+(\mathbf{x}) &= \lim_{z \to x, z \in \Omega} \frac{\partial u^+(\mathbf{z})}{\partial x}, \quad u_y^+(\mathbf{x}) = \lim_{z \to y, z \in \Omega} \frac{\partial u^+(\mathbf{z})}{\partial y} \\ u_x^-(\mathbf{x}) &= \lim_{z \to x, z \in \Omega^c} \frac{\partial u^-(\mathbf{z})}{\partial x}, \quad u_y^-(\mathbf{x}) = \lim_{z \to y, z \in \Omega^c} \frac{\partial u^-(\mathbf{z})}{\partial y} \\ \partial_\mathbf{n} u^+(\mathbf{x}) &= u_x^+(\mathbf{x}) n_x + u_y^+(\mathbf{x}) n_y, \quad \partial_\mathbf{n} u^-(\mathbf{x}) = u_x^-(\mathbf{x}) n_x + u_y^-(\mathbf{x}) n_y \end{aligned}$$

With the denotation above, the jump of variable $u$ and its derivatives at point $\mathbf{x}$ on the boundary $\Gamma$ are defined below

$$\begin{aligned} [u(\mathbf{x})] &\equiv u^+(\mathbf{x}) - u^-(\mathbf{x}) \\ [\partial_\mathbf{n} u(\mathbf{x})] &\equiv \partial_\mathbf{n} u^+(\mathbf{x}) - \partial_\mathbf{n} u^-(\mathbf{x}) \end{aligned}$$

- The volume integral

For a given function $f(\mathbf{x})$ defined on $\Omega$, the volume integral $u(\mathbf{x}) = (Yf)(\mathbf{x})$ defined in (3.16) can be presented as

$$u(\mathbf{x}) = (Yf)(\mathbf{x}) := \int_\Omega G(\mathbf{y}, \mathbf{x}) f(\mathbf{y}) d\mathbf{y} \text{ for } \mathbf{x} \in \mathbb{R}^2 \quad (3.23)$$

The volume integral $u(\mathbf{x})$ can be considered as a solution to the following interface problem

$$\begin{aligned} \Delta u(\mathbf{x}) - \kappa u(\mathbf{x}) &= \tilde{f}(\mathbf{x}) &&\text{for } \mathbf{x} \in \Omega \cup \Omega^c \\ [u(\mathbf{x})] &= 0 &&\text{for } \mathbf{x} \in \Gamma \\ [\partial_\mathbf{n} u(\mathbf{x})] &= 0 &&\text{for } \mathbf{x} \in \Gamma \\ u(\mathbf{x}) &= 0 &&\text{for } \mathbf{x} \in \partial \mathcal{B} \end{aligned} \quad (3.24)$$

The source term $\tilde{f}(\mathbf{x})$ can be written as

$$\tilde{f}(\mathbf{x}) := \begin{cases} f(\mathbf{x}) & \text{for } \mathbf{x} \in \Omega \\ 0 & \text{for } \mathbf{x} \in \Omega^c \equiv \mathcal{B} \setminus \overline{\Omega} \end{cases} \quad (3.25)$$

- The double layer boundary integral

For a given density function $\varphi$ defined on $\Gamma$, the double layer boundary integral $u(\mathbf{x}) = (W\varphi)(\mathbf{x})$ defined in (3.14) can be expressed as

$$u(\mathbf{x}) := (W\varphi)(\mathbf{x}) = \int_\Gamma \frac{\partial G(\mathbf{y},\mathbf{x})}{\partial \mathbf{n_y}} \varphi(\mathbf{y}) ds_\mathbf{y} \text{ for } \mathbf{x} \in \Omega \cup \Omega^c \quad (3.26)$$

The double layer boundary integral $u(\mathbf{x})$ can be considered as a solution to the following interface problem

$$\begin{aligned} \Delta u(\mathbf{x}) - \kappa u(\mathbf{x}) &= 0 &&\text{for } \mathbf{x} \in \Omega \cup \Omega^c \\ [u(\mathbf{x})] &= \varphi(\mathbf{x}) &&\text{for } \mathbf{x} \in \Gamma \\ [\partial_\mathbf{n} u(\mathbf{x})] &= 0 &&\text{for } \mathbf{x} \in \Gamma \\ u(\mathbf{x}) &= 0 &&\text{for } \mathbf{x} \in \partial \mathcal{B} \end{aligned} \quad (3.27)$$

- The single layer boundary integral

For a given density function $\psi$ defined on $\Gamma$, the single layer boundary integral $u(\mathbf{x}) = (V\psi)(\mathbf{x})$ defined in (3.15) can be expressed as

$$u(\mathbf{x}) := (V\psi)(\mathbf{x}) = \int_\Gamma G(\mathbf{y},\mathbf{x})\psi(\mathbf{y}) ds_\mathbf{y} \text{ for } \mathbf{x} \in \Omega \cup \Omega^c \quad (3.28)$$

The single layer boundary integral $u(\mathbf{x})$ can be considered as a solution to the following interface problem

$$\begin{aligned} \Delta u(\mathbf{x}) - \kappa u(\mathbf{x}) &= 0 &&\text{for } \mathbf{x} \in \Omega \cup \Omega^c \\ [u(\mathbf{x})] &= 0 &&\text{for } \mathbf{x} \in \Gamma \\ [\partial_\mathbf{n} u(\mathbf{x})] &= \psi(\mathbf{x}) &&\text{for } \mathbf{x} \in \Gamma \\ u(\mathbf{x}) &= 0 &&\text{for } \mathbf{x} \in \partial \mathcal{B} \end{aligned} \quad (3.29)$$

### 3.2.3 Reinterpretation of the interface problem

In this section, we rewrite the equivalent interface problems (3.24), (3.27), (3.29) for integrals (3.23), (3.26), (3.28) of different types in the unified form, which can be expressed as

$$\begin{aligned}
\Delta u(\mathbf{x}) - \kappa u(\mathbf{x}) &= F(\mathbf{x}), & \mathbf{x} \text{ in } \Omega \cup \Omega^c \\
[u(\mathbf{x})] &= \Phi(\mathbf{x}), & \mathbf{x} \text{ on } \Gamma \\
[u_\mathbf{n}(\mathbf{x})] &= \Psi(\mathbf{x}), & \mathbf{x} \text{ on } \Gamma \\
u(\mathbf{x}) &= 0, & \mathbf{x} \text{ on } \partial\mathcal{B}
\end{aligned} \quad (3.30)$$

| Integral | $\mathcal{F}$ | $\Phi$ | $\Psi$ |
|---|---|---|---|
| $-V\psi$ | $\mathcal{F} \equiv 0$ | $\Phi \equiv 0$ | $\Psi = \psi$ |
| $W\varphi$ | $\mathcal{F} \equiv 0$ | $\Phi = \varphi$ | $\Psi \equiv 0$ |
| $Yf$ | $\mathcal{F} = \tilde{f}$ | $\Phi \equiv 0$ | $\Psi \equiv 0$ |

Obviously, $u(\mathbf{x})$ to the interface problem above can be regarded as the sum of the volume integral, the double layer boundary integral and the single layer boundary integral:

$$u(\mathbf{x}) = (Yf)(\mathbf{x}) + (W\varphi)(\mathbf{x}) - (V\psi)(\mathbf{x}).$$

Once given $f$ and the initial value of $\varphi$ or $\psi$, $u(\mathbf{x})$ defined at grid nodes can be calculated by the interface problem. Therefore, the values specified at boundary $\Gamma$ need to be extracted in each iteration step (3.21), (3.22), which will be discussed in the following subsection. Since the coefficients of the interface problem are constant, it can be solved by fast elliptic solvers such as FFT-based solver or the geometry multigrid method efficiently.

### 3.2.4 Correction

For simplicity, assume the regular domain $\mathcal{B} = [-L, L]^2$ is a square domain. Suppose $M$ is an integer greater than 1, and $h = \frac{2L}{M}$ is the step size of the space. For $i, j = 0, 1, \cdots, M$, assume $x_i = -L + ih$ and $y_j = -L + jh$. The box $\mathcal{B}$ is partitioned into a uniform grid with nodes $\mathbf{p}_{i,j} = (x_i, y_j)$. In order to apply the FFT-based elliptic solver, the five-point finite difference method is applied to discrete the modified Helmholtz equation (3.13) subject to Dirichlet boundary (3.18).

$$\frac{u(x_{i+1}, y_j) + u(x_{i-1}, y_j) + u(x_i, y_{j+1}) + u(x_i, y_{j-1}) - 4u(x_i, y_j)}{h^2} - \kappa u(x_i, y_j)$$
$$= f(x_i, y_j) \quad (3.31)$$

For $i, j = 1, 2, \cdots, M - 1$. Here, $u_{i,j}$ is the finite difference approximation of $u(\mathbf{p}_{i,j})$. The finite difference equation involves a five-point stencil at each grid node $\mathbf{p}_{i,j}$. Obviously, if the righthand side $f$ is sufficiently smooth, then the solution to the finite difference equation (3.31) has second-order accuracy. However, the irregular domain leads to discontinuous $f((\tilde{x}))$, $u((x))$ and $\partial_\mathbf{n} u(\mathbf{x})$ in (3.25), (3.27), (3.29) respectively. As a result, the local truncation error is unable to maintain the second order.

In detail, the grid nodes can be classified into regular nodes and irregular nodes. The grid node whose neighbour grid nodes are located at either side of $\Gamma$ is defined as an irregular node, while others are classified as regular nodes. Due to the discontinuities of the solution and the normal derivative across the interface $\Gamma$, the finite difference approximation at irregular grid nodes has more significant local truncation errors than those at regular grid nodes. The exact estimate is as follows

$$\Delta u(\mathbf{p}_{i,j}) - \kappa u(\mathbf{p}_{i,j}) - f(\mathbf{p}_{i,j}) = \begin{cases} \mathcal{O}(h^2) & \text{if } \mathbf{p}_{i,j} \text{ is a regular node} \\ \mathcal{O}(h^{-2}) & \text{if } \mathbf{p}_{i,j} \text{ is an irregular node} \end{cases}$$

Here, $u(\mathbf{x})$ is the solution to the elliptic PDE (3.30), constrained by two interface conditions on $\Gamma$. The solution to the discrete interface problem is inaccurate since local truncation errors at irregular grid nodes are too large to be ignored.

Due to the existence of the interface in rectangle $\mathcal{B}$, the finite difference equations have large local truncation errors at irregular grid nodes near the interface and thus have to be appropriately modified. The local truncation error at irregular grid recovers to first-order accuracy is sufficient to guarantee that the global solution error of the problem has second-order accuracy[71]. The procedure of correction will be presented in the horizontal and vertical directions as follows.

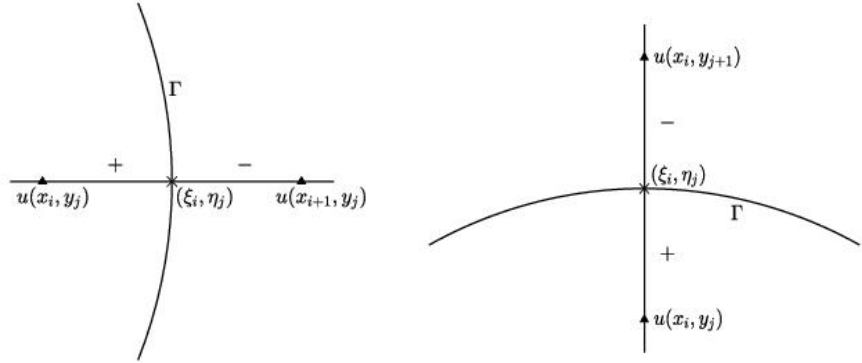

Figure 2: Two cases of the intersecting node in the structure grid. (a) means node lies in horizontal line while (b) implies node lies in vertical line.

- Horizontal direction

Assuming the $\Gamma$ has an intersecting node $(\xi_i, \eta_j)$ located between an exterior node $\mathbf{p}_{i,j}$ and an interior node $\mathbf{p}_{i+1,j}$. The local truncation error of the horizontal direction can be expressed as

$$E_{h,x}(x_i, y_j) \equiv \frac{u^-(x_{i+1}, y_j) + u^+(x_{i-1}, y_j) - 2u^+(x_i, y_j)}{h^2} - u^+_{xx}(x_i, y_j)$$

According to Taylor's expansion, the horizontal part $E_{h,x}(x_i, y_j)$ of the local truncation error can be approximated by the following quantity

$$C^+_{h,x}(x_i, y_j) \equiv -\frac{1}{h^2}\left\{[u] + [u_x](x_{i+1} - \xi_i) + \frac{1}{2}[u_{xx}](x_{i+1} - \xi_i)^2\right\} + H.O.T$$

which introduces a first-order error, since

$$E_{h,x}(x_i, y_j) = C^+_h(x_i, y_j) + \mathcal{O}(h)$$

The correction term $C^+_h(x_i, y_j)$ depends on jumps $[u]$, $[u_x]$ and $[u_{xx}]$. As a matter of fact, the jumps of partial derivatives of the function $u$ can be deduced from the interface problem (3.30). Details of the calculation of jumps are presented in the next subsection. Similarly, $C^-_{h,x}(x_i, y_j)$ denotes the correction of the node $\mathbf{p}_{i,j} \in \Omega^c$ locating near the interface $\Gamma$

$$C^-_{h,x}(x_{i+1}, y_j) \equiv \frac{1}{h^2}\left\{[u] + [u_x](x_i - \xi_i) + \frac{1}{2}[u_{xx}](x_i - \xi_i)^2\right\} + H.O.T$$

- Vertical direction

Assume the $\Gamma$ has an intersecting node $(\xi_i, \eta_j)$ located between an exterior node $\mathbf{p}_{i,j}$ and an interior node $\mathbf{p}_{i,j+1}$. The local truncation error of the vertical direction can be expressed as follows

$$E_{h,y}(x_i, y_j) \equiv \frac{u^-(x_i, y_{j+1}) + u^+(x_i, y_{j-1}) - 2u^+(x_i, y_j)}{h^2} - u^+_{yy}(x_i, y_j)$$

According to Taylor's expansion, the vertical part $E_{h,y}(x_i, y_j)$ of the local truncation error can be approximated by the following quantity

$$C^+_{h,y}(x_i, y_j) \equiv -\frac{1}{h^2}\left\{[u] + [u_y](y_{j+1} - \eta_j) + \frac{1}{2}[u_{yy}](y_{j+1} - \eta_j)^2\right\} + H.O.T$$

which introduces a first-order error, since

$$E_{h,y}(x_i, y_j) = C^+_{h,y}(x_i, y_j) + \mathcal{O}(h)$$

Same to $E_{h,x}(x_i, y_j)$, $E_{h,y}(x_i, y_j)$ also depends on jumps of $[u]$ and $[u_y]$. $C^-_{h,y}(x_i, y_j)$ denotes the correction of the node $\mathbf{p}(x_i, y_j) \in \Omega^c$ locating near the interface $\Gamma$

$$C^-_{h,y}(x_i, y_j) \equiv \frac{1}{h^2}\left\{[u] + [u_y](y_j - \eta_j) + \frac{1}{2}[u_{yy}](y_j - \eta_j)^2\right\} + H.O.T$$

### 3.2.5 Calculation of jumps

This section introduces the calculation of the jumps of the function $u(\mathbf{x})$ and its partial derivatives across the interface, which are required not only for the correction of the finite difference equation, but also for the extraction of boundary data from the discrete finite difference decomposition.

Assume $u(\mathbf{x})$ be a piecewise smooth function defined in $\mathcal{B}$, whose normal derivative $\partial_\mathbf{n} u(\mathbf{x})$ and $u(\mathbf{x})$ itself may be discontinuous across the interface. Computation of the jumps for $u$ on $\Gamma$ can be reiterated below

$$\begin{aligned} \Delta u - \kappa u &= \tilde{f} &\text{in } \mathcal{B} \setminus \Gamma &\quad (3.32)\\ [u] &= \varphi &\text{on } \Gamma &\quad (3.33)\\ [u_\mathbf{n}] &= \psi &\text{on } \Gamma &\quad (3.34) \end{aligned}$$

Let $s$ be the arc length parameter of the curve $\Gamma$, $\mathbf{t} = (\tau_1(s), \tau_2(s))^T$ be a unit tangent vector at a point on the interface. Thus, $\mathbf{n} = (\tau_2(s), -\tau_1(s))^T = \mathbf{t}^\perp$ be the unit outward normal vector on $\Gamma$. Suppose $\varphi = \varphi(s)$ and $\psi = \psi(s)$ are two at least twice differentiable functions defined on $\Gamma$.

Differentiating (3.33) with respect to $s$, together with (3.34):

$$\begin{cases} \tau_1[u_x] + \tau_2[u_y] = \varphi_s \\ \tau_2[u_x] - \tau_1[u_y] = \psi_s. \end{cases} \quad (3.35)$$

First-order partial derivatives of $u$ on $\Gamma$ can be obtained by solving this simple linear system. Continue to differentiate the (3.35) with respect to s, together with identity (3.32), we get

$$\begin{cases} \tau_1^2[u_{xx}] + 2\tau_1\tau_2[u_{xy}] + \tau_2^2[u_{yy}] = \varphi_{ss} - (\tau_1'[u_x] + \tau_2'[u_y]), \\ \tau_1\tau_2[u_{xx}] + (\tau_2^2 - \tau_1^2)[u_{xy}] - \tau_1\tau_2[u_{yy}] = \psi_s - \tau_2'[u_x] + \tau_1'[u_y], \\ [u_{xx}] + [u_{yy}] = f + \kappa[u]. \end{cases}$$

Solving this three-by-three linear system by any direct method or iterative method yields the jumps of second derivatives of $u$ on $\Gamma$.

### 3.2.6 Boundary extraction

Assume $U(\mathbf{x})$ is the numerical solution of $u(\mathbf{x})$. Thus, $U^+, U^-$ denotes the approximate solution of $u^+, u^-$, respectively. Fig. 3 shows the six-point stencil of interpolated point $\mathbf{z}_i$ located in the shadow region of a grid cell. The stencil consists of 6 points and can be obtained by rotation or reflection transformation if $\mathbf{z}_i$ in another grid cell.

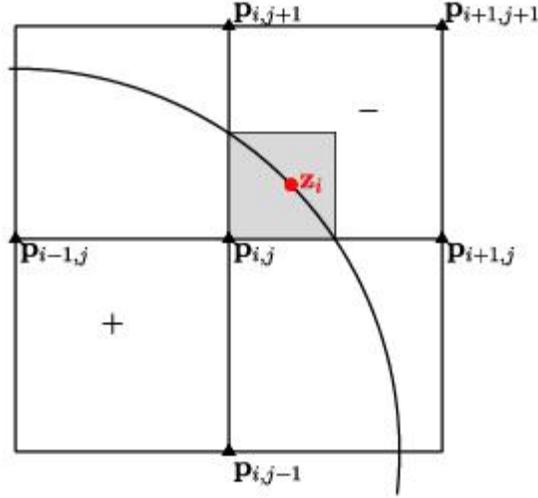

Figure 3: Interpolate stencil

As shown in Fig. 3, $\mathbf{p}_{i,j}$ are the grid points for extracting value of a point $\mathbf{z}_i \in \Gamma$. Following the notations, let $(\xi_i, \eta_j) = \mathbf{p}_{i,j} - \mathbf{z}_i$, thanks to Taylor expansion of $U(\mathbf{z}_i)$ at boundary node $\mathbf{x}$:

$$\begin{aligned}
U(\mathbf{p}_{i,j}) &= U^+(\mathbf{z}_i) + U_x^+(\mathbf{z}_i)\xi_i + U_y^+(\mathbf{z}_i)\eta_j + \frac{1}{2}U_{xx}^+(\mathbf{z}_i)\xi_i^2 \\
&\quad + U_{xy}^+(\mathbf{z}_i)\xi_i\eta_j + \frac{1}{2}U_{yy}^+(\mathbf{z}_i)\eta_j^2 + O(|\mathbf{p}_{i,j} - \mathbf{z}_i|^3), \quad \mathbf{p}_{i,j} \in \Omega \\
U(\mathbf{p}_{i,j}) &= U^-(\mathbf{z}_i) + U_x^-(\mathbf{z}_i)\xi_i + U_y^-(\mathbf{z}_i)\eta_j + \frac{1}{2}U_{xx}^-(\mathbf{z}_i)\xi_i^2 \\
&\quad + U_{xy}^-(\mathbf{z}_i)\xi_i\eta_j + \frac{1}{2}U_{yy}^-(\mathbf{z}_i)\eta_j^2 + O(|\mathbf{p}_{i,j} - \mathbf{z}_i|^3), \quad \mathbf{p}_{i,j} \in \Omega^c
\end{aligned} \quad (3.36)$$

Recall that the boundary data $U^+(\mathbf{x})$, $U_x^+(\mathbf{x})$ and $U_y^+(\mathbf{x})$ are expected. If $\mathbf{p}_{i,j} \in \Omega$, the jump calculated above to correct $U(\mathbf{p}_{i,j})$, $\mathbf{p}_{i,j} \in \Omega^c$

$$J(\mathbf{p}_{i,j}) = [U] + [U_x]\xi_i + [U_y]\eta_j + \frac{1}{2}[U_{xx}]\xi_i^2 + [U_{xy}]\xi_i\eta_j + \frac{1}{2}[U_{yy}]\eta_j^2, \quad \mathbf{p}_{i,j} \in \Omega^c$$

The expansion (3.36) can be reinterpreted as follows:

$$U(\mathbf{p}_{i,j}) + J(\mathbf{p}_{i,j})$$
$$= U^+(\mathbf{z}_i) + U_x^+(\mathbf{z}_i)\xi_i + U_y^+(\mathbf{z}_i)\eta_j + \frac{1}{2}U_{xx}^+(\mathbf{z}_i)\xi_i^2 + U_{xy}^+(\mathbf{z}_i)\xi_i\eta_j + \frac{1}{2}U_{yy}^+(\mathbf{z}_i)\eta_j^2, \quad \mathbf{p}_{i,j} \in \Omega^c$$

Thus, a 6x6 linear system is formed and is always invertible. The new variables can be solved by a simple direct method such as QR decomposition or LU decomposition. $U^+$, $U_x^+$ and $U_y^+$ are recovered by straightforward manipulation of $U^+$, $U_x^+$ and $U_y^+$, respectively.

## 3.3 Algorithm summary

In this section, we summarize the algorithm of the heat, wave and Schrödinger equations as follows:

- Procedure 1: Initialize the structured grid

    a. use quasi-uniformly spaced points $z_i$ to discretize the interface $\Gamma$;

    b. partition $\mathcal{B}$ into a uniform Cartesian grid $\mathcal{T}_h$;

    c. identify the regular and irregular nodes of the Cartesian grid;

    d. find intersecting points located between $\Gamma$ and Cartesian grid line.

- Procedure 2: Evaluate boundary data on control nodes

    a. compute jumps of partial derivatives at control nodes;

    b. solve interface problems 3.30.

    c. extract boundary data $u^+(\mathbf{x})$ or $\partial_\mathbf{n} u(\mathbf{x})$ for Dirichlet or Neumann BVP respectively;

- Procedure 3: Time evolution

    a. initialize the structure grid by Procedure 1 and calculate $F^1$ by initial boundary condition;

    b. choose an arbitrarily initial guess $\varphi^{(0)}$ or $\psi^{(0)}$ to start iteration (3.21) or (3.22) for Dirichlet or Neumann BVP respectively;

    c. calculate the $u^{n+1}(\mathbf{x})$ or $\partial_\mathbf{n} u^{n+1}(\mathbf{x})$ and then update $\varphi^{(k)}$ or $\psi^{(k)}$ by Procedure 2; evaluate right hand side $F^{n+1}$ by former information or initial boundary conditions;

    d. Richardson iteration: repeat previous steps 2 to 4 until $\varphi^{(k)}$ or $\psi^{(k)}$ converges in a specific norm;

## 3.4 GPU accelerated kernel-free boundary integral method

### 3.4.1 Test platform

As in the case of a typical CPU, a GPU also has an Arithmetic logic unit, control unit, cache, and DRAM. A modern NVIDIA GPU contains tens of multiprocessors. A multiprocessor consists of 8 scalar processors (SPs), each capable of executing an independent thread. Threads are grouped together as a grid of thread blocks, such that each MP executes one or more thread blocks and each of its SP runs one or more threads within a block. Our experimental platform is GTX 1080 Ti based on NVIDIA CUDA model, which is also the name of the software for programming this

architecture.GTX 1080 Ti GPU contains in total 3584 cores organized in 88 streaming multiprocessors (MPs). It provides 11GB of device memory with a memory bandwidth of 484GB/s, accessible by all its cores and the CPU through the Intel Xeon E5-2650 v4 with six cores.

Groups of 32 threads (warps) within a block are scheduled with the same instruction in a single cycle of clock time. Furthermore, threads within a block can cooperate by sharing data through some shared memory and must synchronize their execution to coordinate the memory accesses. In contrast, there is no synchronization between thread blocks in a grid of blocks. Hence, the kernels only work with the data in the GPU device's memory, and their final results must be transmitted to the CPU. However, the main bottleneck in the CUDA architecture lies in the data transfer between the host (CPU) and the device (GPU). Any transfer between CPU and GPU may reduce the time execution performance. To summarize, it is necessary to limit the communication between the CPU and GPU.

### 3.4.2 GPU-accelerated KFBI method

In this part, we focus on parallelizing the algorithm in Section 3.3 on the GPU platform. For procedure 1, the domain boundary $\Gamma$ and the structured grid $\mathcal{T}_h$ are set up on the host for the first and second points. Since GPU kernel functions only work on data stored in GPU's video memory, after point 2, the information of grid and boundary is transferred from host memories to the corresponding GPU device memories. In point 3, each grid node is assigned a GPU thread, and the regular and irregular nodes are identified by their related threads. Meanwhile, the number of irregular nodes is counted by GPU threads in parallel. In point 4, the irregular grid nodes are treated as a parallel unit, and the intersection information of each irregular node is calculated by its corresponding thread.

In procedure 2, there are two approaches to parallelize point 1. The first approach involves treating the intersecting point as the parallel unit and utilizing the atomic function *atomicAdd()* in CUDA to correct the irregular points. The second approach, decoupled strategy, involves storing the intersection information of each irregular grid point and assigning a GPU thread to compute the jumps at the corresponding intersections and correct them. For this study, the latter approach is chosen. In point 2, the modified linear system can be solved by following three sub-steps. Firstly, the fast Fourier transform is applied to the data, and the Fourier coefficients are calculated using the CUFFT library. Then, the Fourier coefficients are divided by the eigenvalues. Finally, the inverse fast Fourier transform is applied to the data, and the resulting Fourier coefficients are calculated again using the CUFFT library. To optimize memory access performance, each GPU thread block loads a portion of the grid into its shared memory during these sub-steps. In point 3, each boundary point is assigned to a thread, and the corresponding thread calculates the boundary data.

For procedure 3, we developed a kernel for evaluating the right-hand side in point 4. In this kernel, each GPU thread is assigned to one node. Then, for the Richardson iteration in point 5, each procedure inside the main loop will be executed as described above.

Consequently, we have successfully ported our heat, wave, and Schrödinger equation solvers from C to the GPU platform using CUDA. Within the main loop of time advancement, each kernel is invoked by the host thread and executed on the GPU. The parallel execution takes advantage of multiple threads available on the GPU. In practice, we specify that each block consists of 256 threads in our solvers. Then, the thread block size in each grid is as follows:

$$Blocks = \frac{N + Threads - 1}{Threads},$$

where $N$ is the number of total threads and $Threads$ is the size of thread blocks. The size of $N$ is adjusted adaptively for different kernel functions here. For identifying irregular nodes, solving the modified linear system, and updating the right-hand side $F^n$, $N$ is equivalent to the number of grid nodes. For correction, $N$ is equal to the number of irregular nodes. For interpolating and Richardson iteration in every step, $N$ is the quantity of the boundary nodes.

## 4 Numerical results

To study the numerical accuracy and efficiency of the methods mentioned above, in this section, we present the numerical results for the heat equations, the wave equations, and the Schrödinger equations with Dirichlet and Neumann BVPs in an irregular domain. The number of the quasi-uniformly spaced points on the boundary curve to discretize the BIEs is fixed to be $M$ throughout, where $M$ is the grid size along $x$ or $y$ direction. The bounding box $\mathcal{B}$ embedding the domain $\Omega$ for solving the interface problem is specified as square, whose size as well as the curve parameters will be given respectively in the description of each numerical example.

The following examples give the convergent error of the numerical discretization scheme. The error is defined as $(e_{ij})_{N \times N}$ with $e_{ij} = (\mathbf{u}_h)_{ij} - (\mathbf{u}^*)_{ij}$, where $N$ is the number of interior grid nodes, $\mathbf{u}^*$ is the exact solution, $\mathbf{u}_h$ is the numerical solution with step size $h$. Denote by $\|\mathbf{e}_h\|_\infty$ and $\|\mathbf{e}_h\|_2$ the discrete maximum norm and the scaled discrete $l_2$-norm of $e_{ij}$ respectively, i.e.,

$$\|\mathbf{e}_h\|_\infty = \max_{(x_i, y_j) \in \Omega} |e_{ij}|$$

$$\|\mathbf{e}_h\|_2 = \frac{1}{N} \sqrt{\sum_{(x_i, y_j) \in \Omega} |e_{ij}|^2}$$

To check the algorithm's accuracy, we verify the numerical error in each case with the grid refinement. The tolerance for the Richardson method is fixed to be $1E - 8$. The step size, the number of grid points, and the approximate discrete error in different norms are listed in the corresponding table for each case. Numerical results on the Cartesian grid to the problem are also displayed in the plots for each fixed time.

## 4.1 Numerical results in 2D space

**Example 1**. This example solves the initial boundary value problem of the heat equation $u_t = \Delta u$ on the bullet-shaped domain with the Dirichlet boundary condition. The boundary curve $\partial \Omega$ is fitted with a cubic spline. The initial and boundary conditions are chosen so that the exact solution reads

$$u(x, y, t) = e^{-t}\sin(0.6x + 0.8y) \quad (4.1)$$

The bounding box $\mathcal{B}$ for the interface problem is set to be $\mathcal{B} = [-1.5, 1.5] \times [-1.5, 1.5]$ and the boundary condition is set as $u|_{\partial \mathcal{B}} = 0$. The errors at $T = 1.0$ are summarized in Tab. 1. Numerical results when $T = 0.5, 1.0, 1.5, 2.0$ are shown in Fig. 4.

| domain | grid | $64 \times 64$ | $128 \times 128$ | $256 \times 256$ | $512 \times 512$ |
|---|---|---|---|---|---|
| bullet-shaped | time | 0.25 | 0.125 | 0.0625 | 0.03125 |
| | $\|e_h\|_\infty$ | $1.02E-4$ | $1.47E-5$ | $3.29E-6$ | $8.22E-7$ |
| | $\|e_h\|_2$ | $2.53E-5$ | $6.97E-6$ | $1.52E-6$ | $3.83E-7$ |

*Table 1: Numerical results of Example 1: Dirichlet BVP of the heat equation on the bullet-shaped domain.*

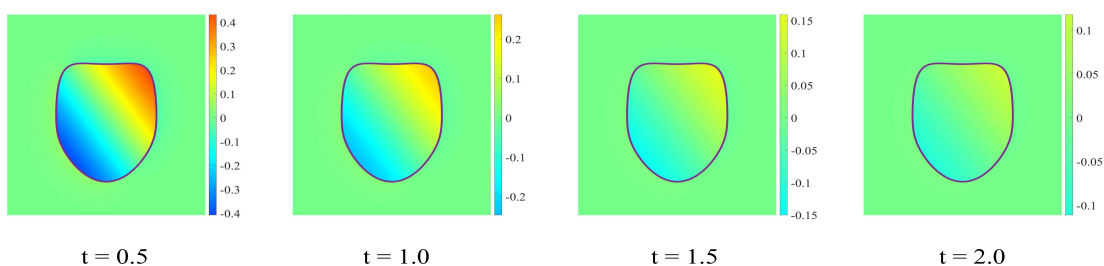

Figure 4: The numerical solutions in Example 1 on the $512 \times 512$ grid.

**Example 2.** This example solves the initial-boundary value problem of the heat equation $u_t = \Delta u$ on the heart-shaped domain with the Neumann boundary condition. The boundary curve $\partial \Omega$ is fitted with a cubic spline. The initial and boundary conditions are chosen so that the exact solution reads

$$u(x, y, t) = e^{-t}\sin(0.6x + 0.8y) \quad (4.2)$$

The bounding box $\mathcal{B}$ for the interface problem is set to be $\mathcal{B} = [-1.5, 1.5] \times [-1.5, 1.5]$ and the boundary condition is set as $u_n|_{\partial \mathcal{B}} = 0$. The errors at $T = 1.0$ are summarized in Tab. 2. Numerical results when $T = 0.5, 1.0, 1.5, 2.0$ are showed in Fig.5.

| domain | grid size | 64 × 64 | 128 × 128 | 256 × 256 | 512 × 512 |
|---|---|---|---|---|---|
| tree-shaped | time step | 0.25 | 0.125 | 0.0625 | 0.03125 |
| | $\|e_h\|_\infty$ | 8.60E − 3 | 2.17E − 3 | 6.21E − 4 | 1.54E − 4 |
| | $\|e_h\|_2$ | 5.00E − 3 | 1.31E − 3 | 3.36E − 4 | 8.24E − 5 |

Table 2: Numerical results of Example 2: Neumann BVP of the heat equation on the heart-shaped domain.

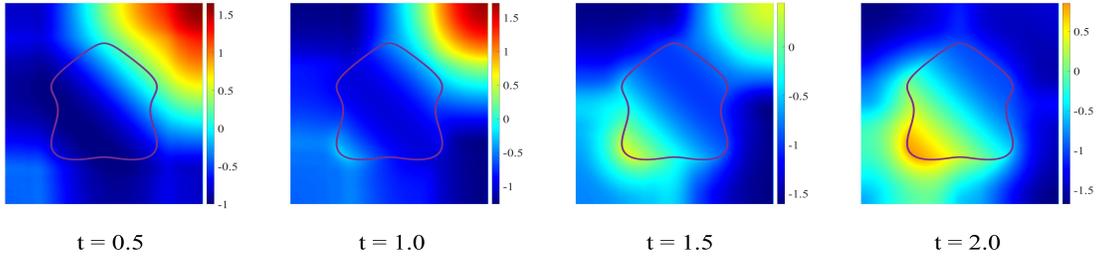

t = 0.5　　　　　t = 1.0　　　　　t = 1.5　　　　　t = 2.0

Figure 5: The numerical solutions in Example 2 on the 512 × 512 grid.

**Example 3**. This example solves the Neumann boundary value problem of the wave equation on the star-shaped domain

$$\begin{cases} x = [(1-c) + c\cos(8\theta)]\cos\theta \\ y = [(1-c) + c\cos(8\theta)]\sin\theta \end{cases} \quad \text{for } \theta \in [0, 2\pi)$$

with the parameter $c = 0.2$. In this example, the model is assumed to satisfy the Neumann boundary condition that $u_n = 0$ on $\partial\Omega$. The initial condition is specified as follows

$$u(x, y, 0) = \frac{1}{1 + \exp(20(\sqrt{x^2 + y^2} - 0.5))}$$
$$u_t(x, y, 0) = 0$$

The bounding box for the interface is set to be $\mathcal{B} = [-1.5, 1.5] \times [-1.5, 1.5]$ and the boundary condition is set as $u_n = 0$ on $\partial\mathcal{B}$. The numerical results depicted in Fig. 6

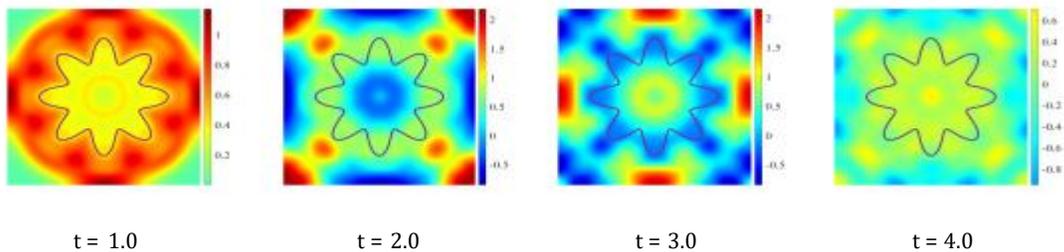

t = 1.0　　　　　t = 2.0　　　　　t = 3.0　　　　　t = 4.0

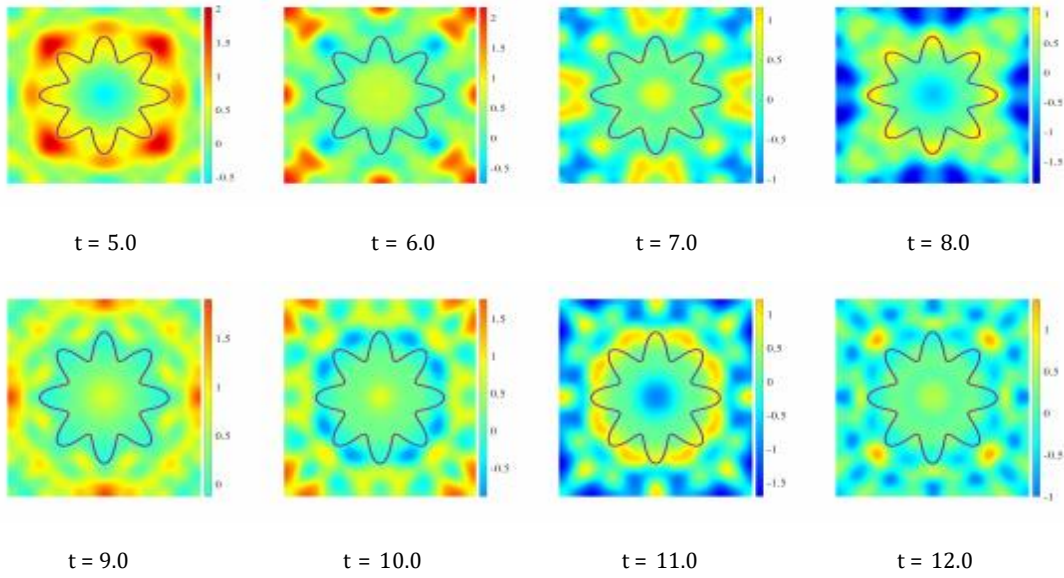

| | t = 5.0 | t = 6.0 | t = 7.0 | t = 8.0 |
| | t = 9.0 | t = 10.0 | t = 11.0 | t = 12.0 |

Figure 6: The numerical solutions in Example 3 on the $512 \times 512$ grid.

illustrate the temporal propagation of waves. Notably, the solutions at each time instance demonstrate an exceptional level of symmetry, providing evidence for the congruence between our algorithm and the underlying physical phenomena.

**Example 4**. According to the papers [72, 73], changing the theta parameter probably changes the stability of the implicit $\theta$−scheme. For the parameter $\theta < 0.25$, the implicit $\theta$−scheme is stable under the CFL condition, while the scheme is unconditionally stable as long as $\theta \geq 0.25$. To verify the stability of the implicit $\theta$−scheme(3.5), the parameter $\theta$ is adjusted in this example. The Dirichlet boundary condition is specified using the exact solution. The boundary curve of the circle domain is given by

$$\begin{cases} x = \cos\theta \\ y = \sin\theta \end{cases} \quad \text{for}\, \theta \in [0, 2\pi)$$

The bounding box for the interface is set to be $\mathcal{B} = [-1.2, 1.2] \times [-1.2, 1.2]$, and the boundary condition is set as $u|_{\partial \mathcal{B}} = 0$. The grid size is set as $256 \times 256$, and the corresponding time step is set as $0.0625$. The errors at $T = 1.0, 10.0, 100.0, 300.0, 500.0$ are shown in Tab.3.

| | T=1 | T=10 | T=100 | T=300 | T=500 |
|---|---|---|---|---|---|
| $\theta = 0.2$ | 4.18E+3 | | | | |
| $\theta = 0.25$ | 6.56E-4 | 1.37E-3 | 2.18E+11 | | |
| $\theta = 0.3$ | 8.49E-4 | 1.80E-3 | 1.25E-2 | 3.05E+7 | |
| $\theta = 0.4$ | 1.23E-3 | 2.73E-3 | 8.59E-4 | 9.60E-1 | 5.67E+4 |
| $\theta = 0.5$ | 1.62E-3 | 9.56E-4 | 1.88E-3 | 1.07E-3 | 2.88E+0 |

Table 3: Numerical results of Example 4: Dirichlet BVP of the wave equation with different $\theta$ in the implicit $\theta$-scheme.

The numerical results summarized in Tab. 3 indicate that when $\theta = 0.2$, the discrete scheme is not satisfied the CFL condition and does not converge at $T = 1$. Convergence over a specific time range can be guaranteed since the implicit $\theta$-scheme is unconditional stability when $\theta \in [0.25, 0.5]$. In addition, from the table, we can conclude that the scheme preserves the accuracy of the numerical algorithm for a longer time as the parameter $\theta$ increases.

**Example 5.** In this example, we consider the Schrödinger equation with $w = 1$, and the potential function [74]

$$v(x, y) = 1 - \cos^2 x \cdot \cos^2 y \quad (4.3)$$

The initial condition

$$u(x, y, 0) = \cos x \cos y$$

Under the condition, the following analytical solution reads

$$u(x, y, t) = \exp(it) \cos x \cos y \quad (4.4)$$

The Dirichlet boundary condition is specified using the exact solution. The boundary curve of the star-shaped domain is defined by

$$\begin{cases} x = 1.5[(1-c) + c\cos(3\theta)]\cos\theta \\ y = 1.5[(1-c) + c\cos(3\theta)]\sin\theta \end{cases} \quad \text{for } \theta \in [0, 2\pi)$$

The bounding box for the interface is set to be $\mathcal{B} = [-\pi, \pi] \times [-\pi, \pi]$ and the boundary condition is set as $u|_{\partial \mathcal{B}} = 0$. The errors at $T = 1.0$ are summarized in Tab. 4. Numerical results when $T = 0, 1.0, 2.0, 3.0$ are shown in Fig. 7.

| domain | grid size | $32 \times 32$ | $64 \times 64$ | $128 \times 128$ | $256 \times 256$ | $512 \times 512$ |
|---|---|---|---|---|---|---|
| | time step | 0.25 | 0.125 | 0.0625 | 0.03125 | 0.015625 |
| star | $\|e_h\|_\infty$ | 1.23E−1 | 3.26E−2 | 8.27E−3 | 2.04E−3 | 4.83E−4 |
| | $\|e_h\|_2$ | 5.73E−2 | 1.55E−2 | 3.91E−3 | 9.62E−4 | 2.30E−4 |

Table 4: Numerical results of Example 5: Dirichlet BVP of the Schrödinger equation on the rotated star-shaped domain.

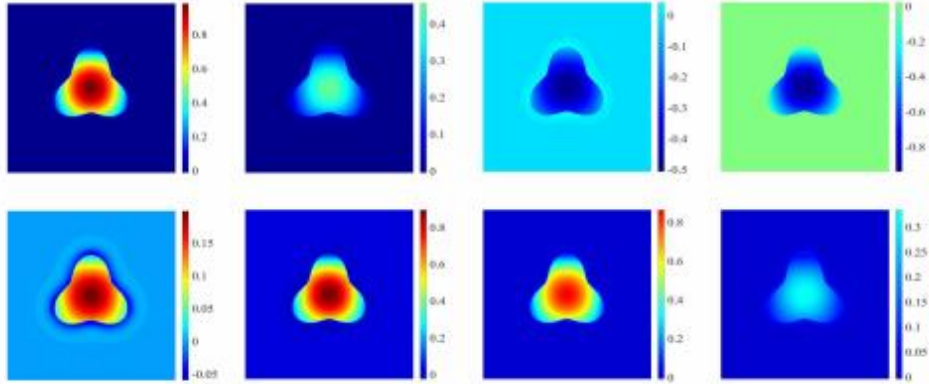

Figure 7: The numerical solutions in Example 5 on the 512 × 512 grid. The first row is the real parts, and the second row is the imaginary parts.

**Example 6.** In this example, we consider the Schrödinger equation with $w = 1$, and the potential function is defined as Eqn (4.3). The initial condition reads

$$u(x, y, 0) = \pi^{-0.25}\exp\left(-0.5(x^2 + y^2)\right)$$

The boundary curve of the star-shaped domain is defined by

$$\begin{cases} x = 1.5[(1 - c) + c\cos(5\theta)]\cos\theta \\ y = 1.5[(1 - c) + c\cos(5\theta)]\sin\theta \end{cases} \text{ for } \theta \in [0, 2\pi)$$

The bounding box $\mathcal{B}$ for the interface problem is set to be $\mathcal{B} = [-\pi, \pi] \times [-\pi, \pi]$ and the boundary condition is set as $u|_{\partial\mathcal{B}} = 0$. The numerical results depicted in Fig. 8 demonstrate the temporal evolution of the real and imaginary parts of the wave function. Additionally, the symmetry of the wave function at each time instance is visually demonstrated in each subfigure.

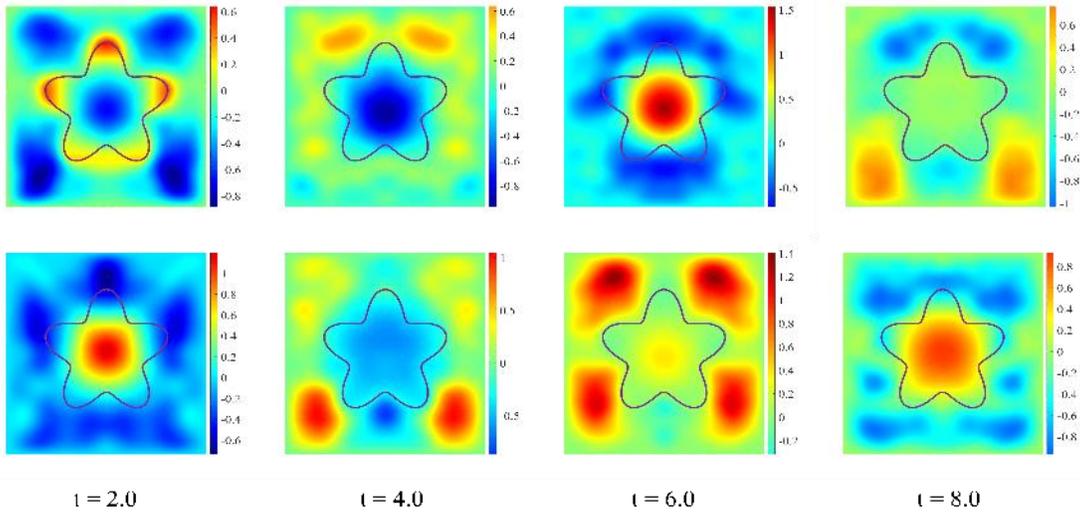

t = 2.0     t = 4.0     t = 6.0     t = 8.0

Figure 8: The numerical solutions in Example 6 on the 512 × 512 grid.

**Example 7.** In this example, we consider the Schrödinger equation with different schemes, including the first-order Godunov splitting method, the second-order Strang splitting method, and the semi-implicit method [75]. The potential function is defined as Eqn (4.3), and the initial and boundary conditions are chosen so that the exact solution satisfies Eqn (4.8). The Dirichlet boundary condition is specified using the exact solution. The boundary curve of the circle domain is given by

$$\begin{cases} x = \cos\theta \\ y = \sin\theta \end{cases} \quad \text{for } \theta \in [0, 2\pi)$$

The bounding box for the interface is set to be $\mathcal{B} = [-1.2, 1.2] \times [-1.2, 1.2]$, and the boundary condition is set as $u|_{\partial \mathcal{B}} = 0$. The grid size is set as $256 \times 256$, and the corresponding time step is set as $0.0625$. The errors at $T = 500.0$ are summarized in Tab. 5,

| Method | T = 500 |
| --- | --- |
| Strang splitting | $1.39E-3$ |
| Semi-implicit | $1.82E-3$ |
| Godunov splitting | $2.73E-2$ |

Table 5: Numerical results of Example 11: Dirichlet BVP of the wave equation with different schemes when $T = 500.0$.

The results in Tab. 5 show that when using these three algorithms to simulate the Schrödinger equation, the phenomenon of the numerical explosion will not be generated. Compared with the time-first operator splitting method, the second-order operator splitting scheme and the semiimplicit scheme have higher accuracy, which is consistent with the theoretical expectation.

In summary, the results in Tab. 1, 2, 4 show that the implicit time discrete scheme combined with the KFBI method has second-order accuracy both in time and space for the heat, wave and Schrödinger equations, independent of the computational domain. The corresponding numerical solutions are shown in Fig.4,5,7. In addition, Numerical results for Example 3, and Example 6 are summarized in Fig. 6, 8, which verifies that our method is still applicable to the problems without exact solution.

**Example 8.** This example shows parallel results for the heat, wave and Schrödinger equations with Dirichlet and Neumann boundary conditions. The boundary curve of the star-shaped domain is defined by

$$\begin{cases} x = 1.0[(1-c) + c\cos(3\theta)]\cos\theta \\ y = 1.0[(1-c) + c\cos(3\theta)]\sin\theta \end{cases} \quad \text{for } \theta \in [0, 2\pi)$$

The initial and boundary conditions are chosen so that the exact solution satisfies (4.1), (4.2), (4.8). The bounding box $\mathcal{B}$ for the interface problem is set to be $\mathcal{B} = [-1.2, 1.2] \times [-1.2, 1.2]$. The terminal time is set as $T = 1$.

Speedup is defined as the ratio between the execution of the time loop on a single CPU $t_{CPU}$ and a single GPU $t_{GPU}$. The CPU time was measured using the Intel Xeon E5-2650 v4 processor, while the GPU time was measured using the GTX 1080 Ti. The specific parameters can be found in Section 3.4.1. The execution times of the examples on the CPU and GPU, as well as the speedups when grid size $N = 128 \times 128, 256 \times 256, 512 \times 512, 1024 \times 1024$ are shown in Tab.10. From the table, starting with low or zero benefits from GPU implementation, increasing the data exploits the parallelization advantages.

| equation | boundary | grid size | $128 \times 128$ | $256 \times 256$ | $512 \times 512$ | $1024 \times 1024$ |
|---|---|---|---|---|---|---|
| Heat | Dirichlet | Time steps | 8 | 16 | 32 | 64 |
| | | CPU time | 0.32 s | 1.22 s | 29.20 s | 235.62 s |
| | | GPU time | 0.36 s | 0.57 s | 1.93 s | 6.98 s |
| | | Speedup | 0.8889 | 2.1579 | 15.1295 | 33.7564 |
| Wave | Neumann | CPU time | 0.64 s | 4.58 s | 38.49 s | 328.85 s |
| | | GPU time | 0.23 s | 0.74 s | 2.18 s | 10.83 s |
| | | Speedup | 2.7826 | 6.1892 | 17.6560 | 30.3647 |
| Schrödinger | Dirichlet | CPU time | 1.47 s | 15.18 s | 162.14 s | 1366.12 s |
| | | GPU time | 0.54 s | 0.80 s | 1.69 s | 6.36 s |
| | | Speedup | 2.7222 | 18.9750 | 95.9408 | 214.7798 |

*Table 6:Simulation time of CPU-based and GPU-accelerated KFBI method, as well as the speedup achieved by GPU-accelerated solver over the CPU-based solver.*

**Example 9**. This example solves the initial boundary value problem of the heat equation with a scale factor $cu_t = \Delta u$ with the Dirichlet boundary condition. The boundary curve of the circle domain is given by

$$\begin{cases} x = \cos\theta \\ y = \sin\theta \end{cases} \quad \text{for } \theta \in [0, 2\pi)$$

the initial and boundary conditions are chosen so that the exact solution satisfies Eqn (4.1), the bounding box $\mathcal{B}$ for the interface problem is set to be $\mathcal{B} = [-1.5, 1.5] \times [-1.5, 1.5]$ and the boundary condition is set as $u|_{\partial \mathcal{B}} = 0$. The terminal time is set as $T = 1$.

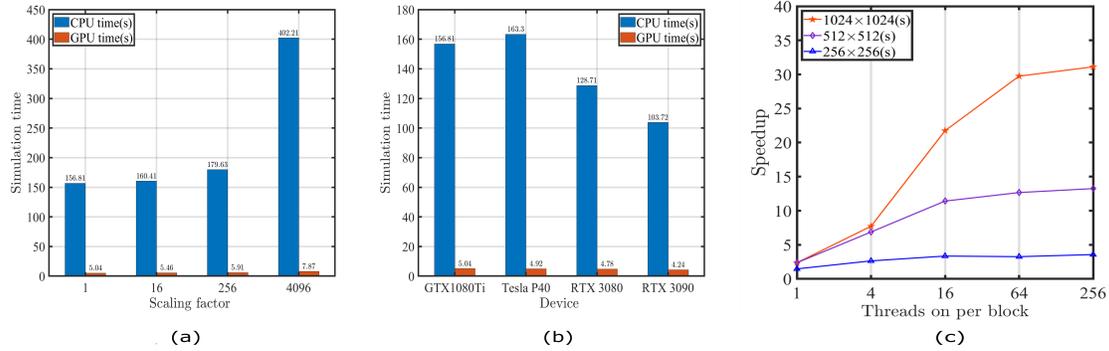

Figure 9: (a)Simulation times of CPU and GPU solvers when scaling factor c = 1,16,256, 4096. Note that the GPU-accelerated solver is 70.19 times faster at the highest scaling factor.(b)Simulation times on four different NVIDIA GPU's with different number of CUDA processors - GTX1080Ti(6 CPU cores and 3584 CUDA cores), Tesla P40(12 CPU cores and 3840 CUDA cores), RTX 3080(12 CPU cores and 8704 CUDA cores) and RTX 3090(16 CPU cores and 10496 CUDA cores).(c)Speedup achieved by our GPU-accelerated solver over the CPUbased solver with varying threads on per block and varying grid sizes.

In Fig. 9(a) we test the effect of parameters $c$ on parallelism efficiency. We select four scale values: 1, 16, 256, and 4096. The average Richardson iteration numbers of all time steps are 13, 13, 14, and 21. In Fig.9(a), we can conclude that the speedup performance is closely related to the scale values. Then we performed scalability analysis of our solver on four different NVIDIA GPUs with varying numbers of CUDA cores: GTX1080Ti, Tesla P40, RTX 3080, and RTX 3090, each with 6, 12, 12, 16 CPU cores, and 3584, 3840, 8704, 10496 CUDA cores. Fig.9(b) shows that as the amount of CUDA cores increases, the computing speed of the GPUaccelerated solver performs faster. Fig.9(c) shows the speedup achieved by our GPU-accelerated solver with different threads per block. As can be seen, the GPU-accelerated solver performs better as the threads on the per block increase under fixed grid size, and the performance of the GPU-accelerated solver scales better as the grid size increases, indicating our method's excellent scalability.

## 4.2 Numerical results in 3D space

The corresponding solvers have been devised for the 3D space of the heat, wave, and Schrödinger equation. The temporal discretization scheme remains consistent with its 2D counterparts. The implicit Laplacian operator is addressed in the spatial direction by adopting the KFBI method. The computational procedures in three-dimensional space differ slightly from those in two-dimensional space, as elaborated in the literature [63].

This section presents numerical results for the Dirichlet boundary value problems of the heat, wave, and Schrödinger equation in 3D space. The computational domain $\Omega$ is enclosed within a bounding box $\mathcal{B}$, which is a simple cube. Detailed information regarding the domain boundary's level set function, the bounding box's size, and the iteration tolerance will be provided for each example. The Numerical results are listed in Tab.7, 8,9, and the color-mapped numerical solution of each example are shown in Fig. 10, 11.

**Example 10**. This example solves the Dirichlet BVP of the heat equation on the four-atom molecular shaped domain $\Omega$ which is given by

$$\Omega = \left\{ \mathbf{x} = (x, y, z) \in \mathbb{R}^3 : c - \sum_{k=1}^{4} \exp\left(-\frac{|\mathbf{x} - \mathbf{x}_k|^2}{r^2}\right) < 0 \right\}$$

with $\mathbf{x}_1 = (\sqrt{3}/3, 0, -\sqrt{6}/12)$, $\mathbf{x}_2 = (-\sqrt{3}/6, 0.5, -\sqrt{6}/12)$, $\mathbf{x}_3 = (-\sqrt{3}/6, -0.5, -\sqrt{6}/12)$, $\mathbf{x}_4 = (0, 0, \sqrt{6}/4)$ and $c = 0.6, r = 0.6$. The bounding box $\mathcal{B}$ for the interface problem is $\mathcal{B} = [-1.2, 1.2] \times [-1.2, 1.2] \times [-1.2, 1.2]$. The tolerance for the Richardson iteration method is fixed to be $1E-8$. The Dirichlet BC is chosen so that the exact solution reads

$$u(x, y, z) = sin(x + y + z)\exp(-3t) \quad (4.5)$$

The errors at $T = 1.0$ are summarized in Tab. 7. Numerical results when $T = 0, 0.5, 1.0$ are shown in Fig. 10.

| grid size | $64 \times 64 \times 64$ | $128 \times 128 \times 128$ | $256 \times 256 \times 256$ | $512 \times 512 \times 512$ |
|---|---|---|---|---|
| time step | 0.25 | 0.125 | 0.0625 | 0.03125 |
| $\|e_h\|_\infty$ | $5.69E-4$ | $9.28E-5$ | $1.91E-5$ | $4.81E-6$ |
| $\|e_h\|_2$ | $2.21E-4$ | $2.83E-5$ | $7.05E-6$ | $1.76E-6$ |

Table 7: Numerical results of Example 10: Dirichlet BVP of the heat equation the four-atom molecular shaped domain.

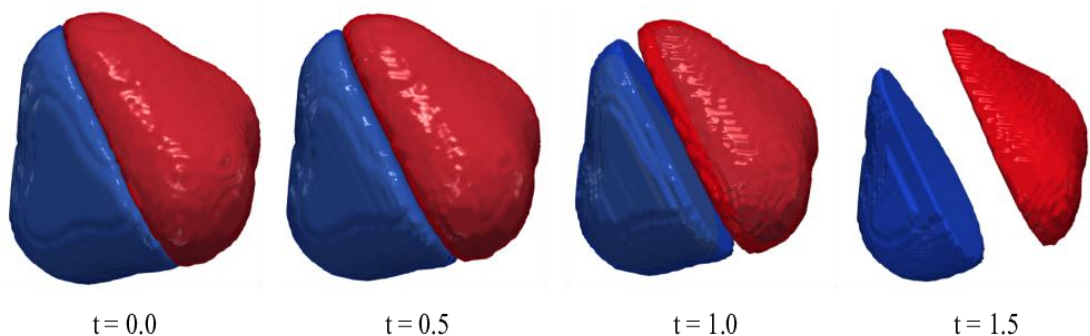

Figure 10: The numerical solutions in Example 10 on the 128 x128x128 grid.

**Example 11**. This example solves the Dirichlet BVP of the wave equation on a tour $\Omega$ which is given by

$$\Omega = \left\{(x, y, z) \in \mathbb{R}^3 : \left(c - \sqrt{x^2 + y^2}\right)^2 + z^2 < a^2\right\}$$

with $a = 0.4$ and $c = 1.0$. The bounding box $\mathcal{B}$ for the interface problem is $\mathcal{B} = [-1.5, 1.5] \times [-1.5, 1.5] \times [-1.5, 1.5]$. The tolerance for the Richardson iteration method is fixed to be $1E - 12$. The Dirichlet BC is chosen so that the exact solution reads

$$u(x, y, z) = sin(x + y + z - \sqrt{3}t) \quad (4.6)$$

The errors at $T = 1.0$ are summarized in Tab. 8. Numerical results when $T = 0, 1.0, 2.0, 3.0$ are shown in Fig. 11.

| grid size | $64 \times 64 \times 64$ | $128 \times 128 \times 128$ | $256 \times 256 \times 256$ | $512 \times 512 \times 512$ |
|---|---|---|---|---|
| time | 0.25 | 0.125 | 0.0625 | 0.03125 |
| $\|e_h\|_\infty$ | $5.90E - 3$ | $1.59E - 3$ | $3.92E - 4$ | $9.88E - 5$ |
| $\|e_h\|_2$ | $2.04E - 3$ | $5.61E - 4$ | $1.48E - 4$ | $3.66E - 5$ |

Table 8: Numerical results of Example 11: Dirichlet BVP of the wave equation on the Torus domain.

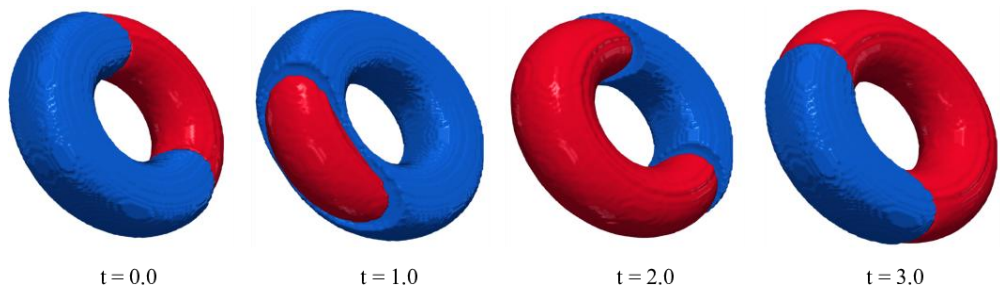

Figure 11: The numerical solutions in Example 11 on the 128 x128x128 grid.

**Example 12**. This example solves the Dirichlet BVP of the Schrödinger equation on an sphere $\Omega$ which is given by

$$\Omega = \left\{(x,y,z) \in R^3 : \frac{x^2}{a^2} + \frac{y^2}{b^2} + \frac{z^2}{c^2} < 1\right\} \quad (4.7)$$

with $a = 1.0, b = 1.0, c = 1.0$. The bounding box $\mathcal{B}$ for the interface problem is $\mathcal{B} = [-1.2,1.2] \times [-1.2,1.2] \times [-1.2,1.2]$. The tolerance for the Richardson iteration method is fixed to be $1E - 12$. The Dirichlet BC is chosen so that the exact solution reads

$$u(x,y,t) = \exp(it)\cos x \cos y \cos z \quad (4.8)$$

The errors at $T = 1.0$ are summarized in Tab.9. Numerical results when $T = 0,1.0,2.0,3.0$ are shown in Fig. 12.

| grid size | $64 \times 64$ | $128 \times 128$ | $256 \times 256$ | $512 \times 512$ |
|---|---|---|---|---|
| time step | 0.25 | 0.125 | 0.0625 | 0.03125 |
| $\|e_h\|_\infty$ | $1.88E-1$ | $1.74E-2$ | $1.81E-3$ | $2.63E-4$ |
| $\|\mathbf{e}_h\|_2$ | $6.91E-2$ | $7.07E-3$ | $7.00E-4$ | $9.83E-5$ |

Table 9: Numerical results of Example 12: Dirichlet BVP of the Schrödinger equation on the sphere domain.

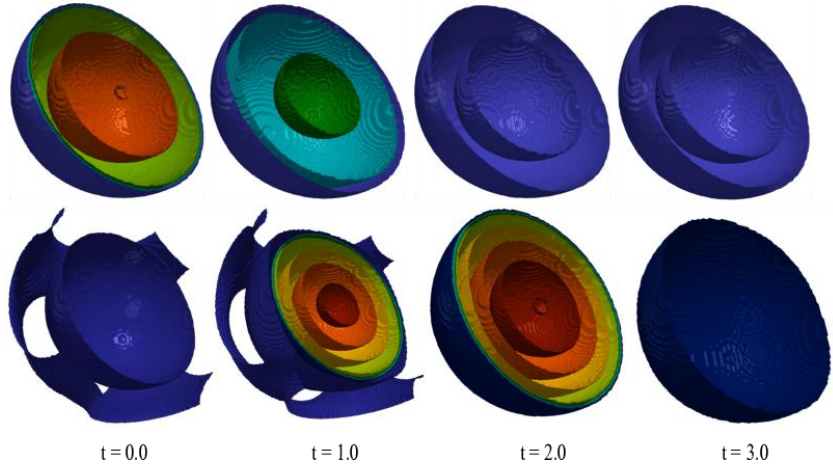

t = 0.0    t = 1.0    t = 2.0    t = 3.0

Figure 12: The numerical solutions in Example 12 on the $128 \times 128 \times 128$ grid. The real part is in the first row and the imaginary part is in the second row.

**Example 13**. This example shows parallel results for the heat, wave and Schrödinger equations in 3D space with Dirichlet boundary conditions. The boundary curve of the sphere domain is defined by (4.7) with $a = 1.0, b = 1.0, c = 1.0$.

For the heat, wave and Schrödinger equations, the initial and boundary conditions are chosen so that the exact solution satisfies (4.5), (4.6), (4.8). The bounding box $\mathcal{B}$ for the interface problem is $\mathcal{B} = [-1.2, 1.2] \times [-1.2, 1.2] \times [-1.2, 1.2]$. The terminal time is set as $T = 1$. The CPU and GPU devices remain consistent with example 13. The execution times of the examples on the CPU and GPU, as well as the speedups when grid size $N = 32 \times 32 \times 32, 64 \times 64 \times 64, 128 \times 128 \times 128, 256 \times 256 \times 256$ are shown in Tab. 10. From the table, starting with low or zero benefits from GPU implementation, increasing the data exploits the parallelization advantages.

| Equation | grid size | $32 \times 32 \times 32$ | $64 \times 64 \times 64$ | $128 \times 128 \times 128$ | $256 \times 256 \times 256$ |
|---|---|---|---|---|---|
| | steps | 2 | 4 | 8 | 16 |
| Heat | CPU time | 1.85 s | 26.36 s | 287.12 s | 4251.62 s |
| | GPU time | 1.04 s | 1.86 s | 3.58 s | 21.51 s |
| | ratio | 1.7788 | 14.1720 | 80.2011 | 197.6578 |
| Wave | CPU time | 1.91 s | 32.17 s | 352.25 s | 4852.81 s |
| | GPU time | 1.02 s | 1.27 s | 4.54 s | 31.38 s |
| | ratio | 1.8725 | 25.3307 | 77.5881 | 154.6466 |
| Schrödinger | CPU time | 3.35 s | 40.09 s | 436.33 s | 6095.76 s |
| | GPU time | 0.72 s | 2.33 s | 7.86 s | 52.51 s |
| | ratio | 4.6527 | 17.2060 | 55.5127 | 116.0876 |

Table 10: DBVP of heat, wave and Schrödinger equations in the sphere domain with radius $r_a = 1.0, r_b = 1.0, r_c = 1.0$, with bounding box $\mathcal{B} = [-1.2, 1.2] \times [-1.2, 1.2] \times [-1.2, 1.2]$ and the tolerance for the Richardson iteration method is fixed to be $1E-12$.

## 5 Discussion

This paper introduces a GPU-accelerated Cartesian grid method with an implicit time discretization scheme for solving the heat, wave, and Schrödinger equations. A second-order implicit discrete scheme is employed for advancing in the temporal dimension. Subsequently, an elliptic boundary value problem is solved within the irregular domain at each time step. A GPU-accelerated KFBI method is developed to circumvent the direct evaluation of boundary integrals. With the parallel KFBI method, the evaluation of a boundary or volume integral is replaced by interpolation of a Cartesian grid-based solution, which satisfies an equivalent interface problem. This interface problem is solved through several procedures, including discretization, correction, solution, and interpolation. The four steps are implemented as kernels using CUDA programming, and each kernel will be executed in parallel in the loop of the Richardson iteration. The numerical solution will be calculated according to the corresponding integral representation formula if the iteration converges. The numerical results represent that this method is an accurate and efficient algorithm applied in smooth irregular areas.

This study solves the discrete boundary integral equations (BIEs) using the parallel Richardson iteration for simplicity. However, alternative Krylov methods, such as the GPU-accelerated GMRES method, can also solve the BIEs. This substitution can significantly reduce the number of iterations, enhancing the efficiency of the method. Moreover, while the CUFFT library is utilized for solving the system on the Cartesian grid in this paper, the Algebraic Multigrid Solver (AmgX) presents an efficient alternative for solving the discrete linear system. Lastly, it should be noted that the GPU-accelerated KFBI method is currently implemented on a single GPU. Future investigations may involve exploring its potential on multiple GPUs and GPU clusters.

The GPU-accelerated KFBI method can be extended into other PDEs such as the Stokes equation, reaction-diffusion equation, Maxwell equation, elasticity equation, etc. Furthermore, when integrated with other algorithms, such as the boundary element or finite element method, the KFBI method may exhibit potential applicability in solving equations within non-smooth domains.

## Acknowledgments

This work is financially supported by the Strategic Priority Research Program of Chinese Academy of Sciences(Grant No. XDA25010405). It is also partially supported by the National Key R&D Program of China, Project Number 2020YFA0712000, the National Natural Science Foundation of China (Grant No. DMS-11771290) and the Science Challenge Project of China (Grant No. TZ2016002).